\documentclass[]{amsart}
\usepackage[]{amsmath, amsthm, amsfonts, rawfonts }

\newcommand {\C} {{\mathbb C}}

\newcommand {\Q} {{\mathbb Q}}
\newcommand {\PP} {{\mathbb P}}

\newcommand {\F} {{\mathcal F}}
\newcommand {\ga} {{\mathcal G \mathcal A}}
\newcommand {\BS} {{\mathbb S}}
\newcommand {\E} {{\tilde E}}

\newtheorem{thm}{Theorem}
\newtheorem{cor}[subsection]{Corollary}
\newtheorem{lemma}[subsection]{Lemma}

\newtheorem{remark}[subsection]{Remark}
\newtheorem{ex}[subsection]{Example}

\begin{document}

\title{
A class of sheaves satisfying Kodaira's vanishing theorem
}
\author{
        Donu Arapura    
}
\address{Department of Mathematics\\
Purdue University\\
West Lafayette, IN 47907\\
U.S.A.}
\thanks{Author partially supported by the NSF}
\email{dvb@math.purdue.edu}
%\date{Feb. 10, 2000} 

\maketitle
%%%%%

This paper contains yet another refinement 
of Kodaira's vanishing theorem; the  result 
in question (corollary \ref{cor:van})
can be viewed as a
Kawamata-Viehweg-Koll\'ar type theorem for vector bundles. 
In order to formulate and prove this  cleanly,
we found it convenient to introduce the class of  {\em geometrically
acyclic} sheaves. It is  the construction and study of
this class (and the related class of {\em geometrically
positive} vector bundles) 
that is the real subject here. 
The name is meant to convey two things: first of 
all, that the sheaves in this class are acyclic in 
the sense that all their higher cohomology groups  vanish. 
Secondly, this class results from  closing
up the subclass of sheaves of ``adjoint type'' under some simple
algebraic and geometric operations.
These closure properties come into play when checking that
a particular sheaf is  geometrically acyclic.
It is our contention that 
many known vanishing theorems can  (and essentially have)
been proved by doing exactly this, and we hope that geometric
acyclicity provides a useful tool in the search for new ones.

We will work throughout over the field of complex numbers $\C$.
The basic constructions and properties of geometrically acyclic sheaves
are worked out in the first  section. The next section
introduces the notion of a geometrically positive vector bundle;
these are bundles for which the naive generalization Kodaira's
vanishing theorem holds. The basic result is that geometrically
positive vector bundles can be obtained by taking 
sufficiently high powers (as described in the appendix)
of vector bundles which are positive in a numerical sense.
The remaining two sections contain various applications, or at any
rate amplifications, of these ideas in algebra and geometry.

My thanks to J. Chipalkatti and J. Lipman for 
helpful conversations, and to  L. Manivel for 
some useful pointers to the literature.

%%%%%%%%%%
\section{Geometric Acyclicity}

\begin{thm}\label{thm:ga}  For each complex projective
variety $X$, there exists a class of coherent sheaves $\ga (X)$ on
$X$, called geometrically acyclic sheaves, satisfying:

\begin{enumerate}
  \item[GA1.] If $Y$ is smooth and projective
 and $f:Y\to X$ is an arbitrary morphism
 then $(f_{*}\omega_Y)\otimes L\in \ga(X)$
for any ample line bundle $L$ on $X$.
\item[GA2.] If 
$$0\to \F_1\to \F_2\to \F_3\to 0$$
is an exact sequence of coherent sheaves on a projective variety
$X$ such that $\F_1,\F_3\in \ga(X)$, then $\F_2\in \ga(X)$.

\item[GA3.] For each $X$, a direct summand of a sheaf in 
$\ga(X)$ is also in $\ga(X)$.
\item[GA4.] If $f: Y\to X$ is an arbitrary
 morphism of projective varieties, then $f_*\F\in \ga(X)$
whenever $\F\in \ga(Y)$.
\item[GA5.] If $\F\in \ga(X)$, 
then $R^ig_*\F = 0$ whenever $i>0$ and 
$g: X\to Z$ is a morphism of projective varieties;
in particular, $\F$  is acyclic.
\end{enumerate}

\end{thm}

\begin{proof}
 Define
$\ga_{0}(X)$ to be the class of sheaves of the type described in
statement GA1.
Inductively define $\ga_{n+1}(X)$ to be the class of sheaves $\F$
such that either $\F$ is an extension of sheaves in $\ga_{n}(X)$,
$\F$ is a direct summand of a sheaf in $\ga_{n}(X)$, or $\F = 
f_{*}\F'$ where $f:Y\to X$ is a morphism of projective varieties
and  $\F'\in \ga_{n}(Y)$. Set
$\ga(X)  = \bigcup_{n}\, \ga_{n}(X)$. The first four statements clearly
hold.

Statement GA5 will be proved by induction. By definition, any sheaf in
$\ga_{0}(X)$ is of the form $\F=(f_{*}\omega_Y)\otimes L$
with $Y$ smooth, $f:Y\to X$ a morphism, and
 $L$ an ample line bundle  on $X$. Given a morphism $g:X\to Z$
with $Z$ projective, choose an ample line bundle $M$ on $Z$.
Then 
$$H^i(X, f_*\omega_Y\otimes L\otimes g^*M^{\otimes N}) = 0$$
for 
$i > 0$ and $N \ge 0$ by Koll\'ar's vanishing theorem \cite{kollar}.
 Therefore by Serre's vanishing
theorem, the Leray spectral sequence and the projection formula,
we obtain
$$H^0(Z, R^ig_*(f_*\omega_Y\otimes L)\otimes M^{\otimes N}) = 0$$
when  $i> 0$ and $N >> 0$ .
As the sheaves above can be assumed to be globally generated (again 
by Serre), we obtain
$R^ig_*(f_*\omega_Y \otimes L) = 0$ for $i > 0$. 

Now suppose that GA5 has been established for sheaves in $\ga_{n}(X)$
(for all $X$, $g$ and $i$). Choose a sheaf $\F\in \ga_{n+1}(X)$
and a morphism $g:X\to Z$. If $\F$ is  a summand or an extension
of sheaves in in $\ga_{n}(X)$, then clearly $R^{i}g_{*}\F=0$ for
$i>0$. Therefore we may assume that $\F = f_{*}\F'$ with $\F'\in 
\ga_{n}(Y)$ with $f:Y\to X$ a morphism. By the induction hypothesis, 
$R^{i}f_{*}\F'=0$ for $i > 0$. Therefore the spectral sequence for the 
composite collapses to yield isomorphisms
$$R^{i}g_{*}f_{*}\F' \cong R^{i}(g\circ f)_{*}\F'$$
But the sheaves on the right also vanish for $i>0$ by induction.

\end{proof}

As should be clear from the proof, $\ga = \bigcup_X\, \ga(X)$
is taken to be the smallest class for which the theorem
holds. We expect that there are in fact  natural extensions
of this class for which properties GA1 to GA5 continue to hold. 
However, this will not be pursued here.

Let us discuss a few basic examples of geometrically acyclic
sheaves.

\begin{ex}
Let $X$ be a smooth projective variety, and let $D\subset X$ be a reduced
divisor with normal crossings. Then $\omega_{X}(D)\otimes L$ is 
geometrically acyclic, provided that $L$ is an ample line bundle.
To see this, write $D = D_{1}+ D'$ where $D_{1}$ is a component.
Then, using the Poincar\'e residue,  $\omega_{X}(D)\otimes L$ can
be expressed as an extension of 
$\omega_{D_{1}}(D')\otimes L$ by $\omega_{X}(D')\otimes L$. Hence
the result follows by induction on the number of components.
    
\end{ex}  

\begin{ex}
Let $X$ be a smooth projective variety. Let $D$ a $\Q$-divisor with 
normal crossing support, i.e. a sum $\sum_{i}a_{i}D_i$ with $a_{i}\in 
\Q$ such that $\sum_{i}D_{i}$ is a normal crossing divisor.
Let $L$ be a line bundle such that $L(D)$ is nef (i.e., 
$c_1(L(D)) :=c_1(L)+ \sum_i a_i[D_i]$ has
a nonnegative intersection number with any curve on $X$)
and big (i.e., $c_1(L(D))^{dimX}>0$).
 Then $\omega_{X}\otimes L(\lceil D\rceil)$ is geometrically acyclic,
where the symbol $\lceil D \rceil$
means that the coefficients should be rounded up to the nearest 
integers. The proof of this is implicit in \cite{kawamata}.
By using embedded resolution of singularities and
a series of covering tricks, he showed
that there exists a generically finite map 
of smooth projective varieties $\pi:X'\to X$, an ample line
bundle $L'$, and  a reduced normal crossing divisor $D'$ on $X'$,
such that $\omega_{X}\otimes L(\lceil D\rceil)$ is a direct
summand of $\pi_{*}(\omega_{X'}\otimes L'(D'))$.
\end{ex}

It will be convenient to extend these definitions to singular
varieties. The extensions will be based on properties which
are known to equivalent in the smooth case (see \cite{kawamata,
mori}).
A line bundle $L$ on a projective variety $X$ is nef
provided that for any ample line bundle $H$ and positive integer
$n$, $L^{\otimes n}\otimes H$ is ample. $L$ is big provided that
for  any ample line bundle $H$, 
$L^{\otimes n}\otimes H^{-1}$ has a nonzero 
global section for some  $n>0$. Alternatively, $L$ is big if and if
$h^0(L^{\otimes n}) = O(n^{dim X})$ (after replacing $L$ by a 
sufficiently high power). In order to check that
$L$ is big it suffices to show that 
$L^{\otimes n}\otimes H^{-1}$ has a nonzero 
global section  for some $n>0$ and a line bundle $H$ which is 
already known to be big.

\begin{lemma}\label{lemma:split} If $\pi:X'\to X$ is a generically finite map of
smooth projective varieties then $\omega_{X}$ is a direct summand
of $\pi_*\omega_{X'}$
\end{lemma}

\begin{proof} 
The result is certainly well known, so we will be quite 
brief. There are two maps: the Grothendieck trace
$\tau:\pi_*\omega_{X'}\to \omega_X$ and a map $p$ in the opposite
direction which corresponds to pullback of forms.
The lemma follows from the identity
 $\tau\circ p = deg \pi$.
\end{proof}

The third example, which will be needed later,
is a slight variation on the previous ones.

\begin{lemma}\label{lemma:ex3}
 Let $f:Y\to X$ be  a surjective morphism of
smooth projective varieties, and $L$ a nef and big line
bundle on $X$, then $f_*\omega_Y\otimes L$ is geometrically
acyclic.
\end{lemma}

\begin{proof}
In the course of the proof, we will construct a large commutative
diagram:
$$
\begin{array}{ccccc}
    \cdots & Y_{i}           & \to&\cdots & Y \\
    \cdots & f_{i}\downarrow &    &\cdots &\downarrow \\
    \cdots & X_{i}           & \to&\cdots & X \\
\end{array}
$$
We will denote the maps $X_i\to X_{i-1}$, $X_{i}\to X$,
$Y_{i}\to Y_{i-1}$ and $Y_{i}\to Y$ by $p_i$, $P_i$,
$\pi_i$ and $\Pi_i$ respectively.
    
Let $H$ be a very ample divisor on $X$. If $n >> 0$,
$L^{\otimes n}(-H)$   has a nonvanishing section, and
let $D$ be the corresponding
effective divisor. Let $X_1\to X$ be a birational map
of smooth varieties 
such that the pullback $D_1$ of $D$ has normal crossings.
Let $Y_1$ be  a
desingularization of $Y\times_{X}X_1$ such that $f^{-1}D_1$
also has normal crossings \cite{hironaka}.

We  can choose  a positive rational linear combination of exceptional
divisors  $E_1$ such that $H_1=P_1^{-1}H-E_1$ is an ample 
$\Q$-divisor. Define
$$H_1' = \epsilon H_1 + (1-\epsilon)P_1^{-1}(H+D)$$
$$ = H_1 + (1-\epsilon)(E_1+ D_1)$$
where $0< \epsilon << 1$ is rational, then $H_1'$ is an ample
$\Q$-divisor. Therefore
$$P_1^{-1}(H+D) = H_1' + \epsilon(E_1+ D_1).$$ 
Therefore, after replacing $n$ by a large multiple, we can
assume that the pullback of the linear system
associated to $L^{\otimes n}$ has an integral normal crossing
divisor  of the
form $H_{1}'' + D_{1}''$ where $H_{1}''$ is smooth, 
ample and transverse to $Y_1\to X_1$, and 
multiplicities of the components of $D_{1}''$ and $f^{-1}D_1''$
are less than $n$. Apply Kawamata's trick \cite[thm 17]{Ka1}, 
to obtain a finite cover $X_2\to X_1$,
branched over $H_1''$ such that the pullback $H_2$ of $H_1''$ has 
multiplicity $n$ and the multiplicities of the components of the 
pullback  $D_2$ of $D_1''$ are unchanged. Let $f_2:Y_2\to X_2$ be the 
fiber product of $Y_1$ with $X_2$.
By our assumptions, $Y_2$ is smooth and 
$f_2^{-1}(H_2+D_2)$ is a normal crossing divisor such that $f_2^{-1}H_2$
has multiplicity $n$ and all other components have smaller multiplicity.

Let $Y_3\to Y_2$ be a desingularization of the $n$-fold cyclic cover  
branched along   $f_2^{-1}(H_2+D_2)$ determined by the line bundle
$f_2^*P_2^*L$ (see  \cite{esnault-viehweg, viehweg}). 
Using the formula for the canonical
sheaf of a cyclic cover [loc. cit.], we see that 
$$\pi_{3*}\omega_{Y_3} = \bigoplus_{i=0}^{n-1}\,
\omega_{Y_2}\otimes f_2^*P_2^*L^{\otimes 
i}(-[\frac{i}{n}f_2^{-1}(H_2+D_2)]) .$$
By our assumptions about multiplicity, the summand on the right
for $i=1$ is just
$\omega_{Y_2}\otimes f_2^{*}P_2^*L(-f_2^{-1}H_2)$.
By lemma \ref{lemma:split}, we see that 
$f_*\omega_Y\otimes L$ is a direct summand of 
$f_*\Pi_{2*}(\omega_{Y_2}\otimes f_2^{*}P_2^*L)$ which is in turn
a direct summand of $P_{2*}((f_2\circ\pi_3)_*\omega_{Y_3}\otimes 
O(H_2))$. Therefore it is geometrically acyclic.

\end{proof}

%%%%%%%%%%%%%%
\section{Geometric Positivity}

Let us say that a vector bundle $E$  is {\em geometrically
semipositive} provided that $\F\otimes E$ is geometrically
acyclic whenever $\F$ is. We will say that a vector bundle
$E$ on $X$ is {\em geometrically positive} provided that it is 
geometrically semipositive, and $f_*\omega_Y\otimes E$ is geometrically
acyclic whenever $f:Y\to X$ is a surjective map from a smooth
projective variety.

\begin{ex} 
A nef line bundle $L$ is 
geometrically semipositive. This follows by induction:
If $\F\in \ga_0(X)$, then $\F\otimes L\in \ga_0(X)\subset\ga(X)$
because $L\otimes H$ is ample whenever $H$ is an ample line bundle.
Suppose that the result is known for all $\F\in \ga_n= 
\bigcup_{Y}\,\ga_n(Y)$, then clearly it holds
for a direct summand  or an extension  of a sheaves in $\ga_n$. If $\F = 
f_*\F'$ with $\F'\in \ga_n$, then $\F\otimes L = 
f_*(\F'\otimes f^*L)\in \ga$. 

This together with lemma \ref{lemma:ex3} shows that a 
 nef and big line bundle is geometrically positive.
\end{ex}

It should be clear that the class of geometrically (semi-)positive
vector bundles is stable under extensions, direct summands,
and tensor products, and therefore also under the Schur-Weyl
powers (see appendix).

\begin{lemma} If $E$ is geometrically semipositive (respectively
positive), then so is $\BS^\lambda(E)$ for any (nonzero) partition
$\lambda$.
\end{lemma}

\begin{proof} $\BS^\lambda(E)$ is a direct summand of a tensor
power of $E$.
\end{proof}

Let us say that a vector bundle $E$ on a projective variety $X$
is nef provided that $O_{\PP(E)}(1)$ is nef.
(We use the convention that $\PP(E)
= \mathbf {Proj}(S^\bullet (E))$.)
It is easy to see that
locally free quotients of nef vector bundles, in particular
globally generated  vector bundles, are nef with this 
definition. Less obvious is the fact that direct sums and
tensor products of
nef bundles are nef. This can deduced easily from the following
lemma and the ampleness of  direct sumas, tensor products and quotients
of ample vector bundles \cite{hartshorne}.

\begin{lemma} Let $H$  be an ample line bundle on $X$. $E$
is nef if and only if $S^n(E)\otimes H$ is ample for all
$n>>0$.
\end{lemma} 

\begin{proof}
Suppose that $E$ satisfies the condition of the lemma,
but that $O_{\PP(E)}(1)$ fails to be nef.
Then there exists a curve $C\subset
  \PP(E)$, with $c_1(O(1))\cdot[C] < 0$.
Since $O(n)\otimes \pi^*H$ is quotient of $\pi^*[S^n(E)\otimes H]$ where
$\pi$ is the projection, 
$$c_1(O(1))\cdot[C] +  \frac{1}{n}\pi^*c_1(H)\cdot[C]\ge 0.$$
This yields a contradiction when $n\to \infty$.

Conversely suppose that $E$ is nef. 
 Nakai-Moishezon's criterion shows that $O(1)\otimes \pi^*H$ 
is ample, therefore $E\otimes H$ is ample.
There  exists a finite
branched cover $p:Y\to X$ such that $p^*H$ possesses an $n$th root
\cite{bloch-gies}, that is a line bundle $L$ such that $L^{\otimes n}
\cong p^*H$.
Therefore  $p^*(S^nE\otimes H)\cong S^n(p^*E\otimes L)$ is ample, 
and this implies
that $S^n E\otimes H$ is ample.
\end{proof}

The property of being nef is also known as numerical
semipositivity. The two notions of semipositivity are related:

\begin{lemma}\label{lemma:gsisnef}
 A geometrically semipositive bundle $E$ on a smooth projective
variety $X$ is nef.
\end{lemma}

\begin{proof} 
Let $O_X(1)$ be a very ample line bundle such that 
$H= \omega_X\otimes O_X(dim\, X + 2)$
is also ample. Then $H^i(X, S^n(E)\otimes H(-1-i)) = 0$ 
for all $i>0$ and $n>0$. In other words, $S^n(E)\otimes H$ 
is $(-1)$-regular \cite[p 100]{mumford}. This implies
that $S^n(E)\otimes H$ is a quotient of a sum of $O_X(1)$'s,
and therefore ample. Thus $E$ is nef by the previous lemma.

\end{proof}

We call a vector bundle $E$ big provided
 that there exists a very  ample line bundle $H$ and
 $n>0$, such that  $S^n(E)\otimes H^{-1}$ is generically generated
by global sections, which means that
$$H^0(S^n(E)\otimes H^{-1})\otimes O_X\to S^n(E)\otimes H^{-1}$$
is surjective over a nonempty open set

 When $E$ is big,  
$S^n(E)\otimes L^{-1}$ is 
generically globally generated for any line bundle $L$ and some
$n>0$, because a negative power of $H$ injects into $L^{-1}$. 
The pullback of an ample bundle under a birational map is both
nef and big.  This notion of bigness coincides with the
previous one for line bundles, and agrees with a more general
one for torsion free sheaves  \cite[pp 292-293]{mori}.
Bigness of $E$ is a stronger condition than
the bigness of  $O_{\PP(E)}(1)$ as the following shows:

\begin{lemma}\label{lemma:big} 
Suppose that  $E$ is a big vector bundle 
over a projective variety $X$. Let $Y$ be projective
variety with a map
$f:Y\to \PP(E)$ which is  generically finite over $f(Y)$
and  such that the composite $Y\to X$ is surjective,
then $f^*O_{\PP(E)}(1)$ is big.
\end{lemma}

\begin{proof}
Denote the projection $\PP(E)\to X$ by $\pi$.  
Let $H$ be an ample line bundle on $\PP(E)$ which is necessarily
of the form $O_{\PP(E)}(a)\otimes \pi^*L$ where
$L$ is a line bundle on $X$. As $f^*H$ is big, 
it is enough to check that
$f^*(O_{\PP(E)}(N)\otimes H^{-1})$ has a nonzero section for
some $N>0$ (see the comments preceding lemma \ref{lemma:split}).

Let $y$ be a general point of $Y$; it lies over a general
point $x\in X$. $\C(x)$ will denote the residue field at $x$.
 Choose a trivialization of  $L$ at $x$.
Then a section of $S^n E_x\otimes L_x^{-1}\otimes \C(x)$
can be identified with
a section of $S^n E_x\otimes \C(x)$, and this determines
a hypersurface in the fiber $\PP(E)_x$. 
Choose $n>0$ so that all elements of $S^n E_x\otimes L_x^{-1}$
extend to global sections, and
choose a section
$\sigma'\in S^n E_x\otimes L_x^{-1}\otimes \C(x)$, 
so that the corresponding hypersurface avoids $f(y)$. 
Then $\sigma'$ lifts to a global section
of $S^n E\otimes L^{-1}$ which  can be identified
with a global section $\sigma$ of $O_{\PP(E)}(n+a)\otimes
H^{-1}$. Pulling $\sigma$ back to $Y$ yields
a nonvanishing section of $f^*(O_{\PP(E)}(n+a)\otimes H^{-1})$.
\end{proof}

\begin{lemma} A locally free quotient of a big vector bundle is big.

\end{lemma}

\begin{proof}
Let $F$ be a locally free quotient of a big vector bundle $E$.
Suppose $H$ is very ample, choose $n>0$ so that  $S^n E\otimes H^{-1}$
is generically generated by its global sections. Then
the restriction of these sections also generates
$S^n F\otimes H^{-1}$ generically.
\end{proof}

\begin{lemma} If $E$ is big, then for any coherent
sheaf $\F$ there exists an integer $n_0$ such that
$S^n(E)\otimes \F$ is generically globally generated
whenever $n \ge n_0$.
\end{lemma}

\begin{proof}
 The proof will be
reduced to a series of observations.
Fix a very ample line bundle $H$.

The tensor product of two or more generically globally generated
coherent sheaves has the same property. Therefore
the set $I(H)$ of integers $n$
for  which  $S^n(E)\otimes H^{-1}$ is generically globally
generated forms a semigroup.

Given a coherent sheaf $\F$,  choose an integer
such that  $\F\otimes H^{\otimes m}$ is globally generated.
Then
$S^n(E)\otimes \F$ is generically
globally generated for any 
$n\in I(H^{\otimes m})\subset I(H)$.

Apply the result of the previous paragraph to obtain
$n\in I(H)$ such that $S^n(E)\otimes E\otimes H^{-1}$
is generically globally generated. This implies that 
${n+1}\in I(H)$. A semigroup containing two
relatively prime integers contains all but finitely
many positive integers, and this concludes the
proof.

\end{proof}

\begin{lemma} The direct sum of two big vector bundles
is big.
\end{lemma}

\begin{proof}
Let $E_1$ and $E_2$  be two big vector bundles
and $H$ a very ample line bundle. Then for each $m\ge0$,
there is an  integer $N_m   > 0$ such that for both  $i$
and all $n \ge N_m$, 
$S^{n}(E_i)\otimes H^{-\otimes m}$ is generically globally
generated. 
Choose $M>0$ so that 
$S^{n}(E_i)\otimes H^{\otimes M}$ is globally
generated for all $n < N_1$. Then  choose $R \ge 2N_{M+1}$.
Then one sees, after grouping terms appropriately, that
$$S^n(E_1)\otimes S^{R-n}(E_2)\otimes H^{\otimes - (1+M)}\otimes
H^{\otimes M}$$
is  generically globally generated for any $n\ge 0$.
Therefore the same holds for $S^R(E_1\oplus E_2) \otimes H^{-1}$. 

\end{proof}

\begin{lemma}
A Schur-Weyl power of a big vector bundle is big.
The tensor product of two big vector bundles is big.
\end{lemma}    

\begin{proof}
Suppose that $E$ is a big vector bundle,
and $H$ an ample line bundle.
For any partition $\lambda$, there exists 
integers
$N_1,\ldots N_r$ such that 
 $min(N_i)\to \infty$ as the weight $|\lambda|\to \infty$ and 
 such that $\BS^\lambda(E)$
is a direct summand of 
$$S^{N_1}(E)\otimes \ldots \otimes S^{N_r}(E)$$
by \cite[5.1]{hartshorne}. This implies that
$\BS^\lambda(E)\otimes H^{-1}$ is
generically globally generated for $|\lambda|>>0$.
Therefore $S^N(\BS^\lambda(E))\otimes H^{-1}$ is generically
globally generated for arbitrary $\lambda$ and $N>>0$,
 since the first factor  can be decomposed into a sum of 
Schur-Weyl powers of large weight.

If $E$ and $F$ are both big, then $E\otimes F$
must be big since it is a direct summand of $S^2(E\oplus F)$.

\end{proof}

The converse to lemma \ref{lemma:gsisnef}
is false
 as we will see (example \ref{ex:tangentb}).
 However, the following may be viewed as a weak converse.
(The notation is explained in the appendix.)

\begin{thm}\label{thm:positiv}
Let $E$ be a vector bundle on a projective variety 
and $\lambda$ a nonzero  partition.
\begin{enumerate}
\item If $E$ is nef, then $\BS_+^\lambda(E)$ is geometrically
semipositive.
\item If $E$ is nef and big, then $\BS_+^\lambda(E)$ is geometrically
positive.
\end{enumerate}
\end{thm}

Before giving the proof, we will need the following lemmas. 
 
\begin{lemma}\label{lemma:canonical}
 Suppose that $ k_1<\ldots <k_m$ is a sequence
of positive integers.
Let $E$ be a rank $k_{m+1}$ vector bundle over a
variety $X$, let $p:F= Flag_{k_1,k_2\ldots k_m}(E)\to X$ be the
bundle of partial flags, and let $\pi_i:F\to Grass_{k_i}(E)$
be the natural projections. If we set $k_0=0$, then
$$\omega_F = p^*\omega_{X}\otimes p^*(det E)^{\otimes k_m}\otimes 
\pi_1^*O_{Grass_{k_1}}(k_0-k_2)\otimes \ldots
\pi_m^*O_{Grass_{k_m}}(k_{m-1}-k_{m+1}) $$
where $O_{Grass_{k_i}}(1)$ are the restrictions of the hyperplane
bundles under the Pl\"ucker embeddings.
\end{lemma}

\begin{proof} See \cite[2.10]{demailly} \end{proof}

\begin{lemma}\label{lemma:nef}
 Let $E, F, X$ be as in the previous lemma, and
$a_1, \ldots a_m$ be positive integers. Then 
$\pi_1^*O_{Grass_{k_1}}(a_1)\otimes\ldots  
\pi_m^*O_{Grass_{k_m}}(a_m)$
is nef (respectively nef and big) if $E$ is nef (respectively
nef and big).
\end{lemma}

\begin{proof}

Let $P=\PP(\wedge^{k_1}E)\times_{X}\ldots\times_{X}\PP(\wedge^{k_m}E) $
and let $i:F\hookrightarrow P$ 
denote the Pl\"ucker embedding.
The restriction of $M(a_1,\ldots, a_m) = O(a_1)\otimes\ldots O(a_m)$
to $F$ is $\pi_1^*O_{Grass_{k_1}}(a_1)\otimes\ldots  
\pi_m^*O_{Grass_{k_m}}(a_m)$.
If $E$ is nef, then so is
$M(a_1,\ldots, a_m)$,
 and hence also $i^*M(a_1,\ldots a_m)$.
 This argument requires only that all $a_j\ge 0$.

Suppose that $E$ is nef and big, then
$E'= \wedge^{k_1}E\oplus\ldots \oplus\wedge^{k_m}E$ is also nef and big
by the preceding lemmas.
Let $j:P\hookrightarrow \PP(E')$ be the Segre embedding.
Then $i^*j^*O_{\PP(E')}(1)= i^*M(1,\ldots 1)
$ is nef and big by lemma \ref{lemma:big}. This together
with the previous paragraph implies that
$$i^*M(a_1, \ldots a_m) = 
i^*M(1,\ldots 1)\otimes i^*M(a_1-1, \ldots a_m-1)$$
is again nef and big.
\end{proof}

\begin{proof}[Proof of theorem]
Let $k_1<\ldots <k_m$ be the list of indices $k$ for which  
$ \lambda_k-\lambda_{k+1}\not=0$, and let 
$a_i = \lambda_{k_i}-\lambda_{k_i+1}$.
Also let $k_{m+1} = rk(E)$ and $k_0=0$. 

Let $F=Flag_{k_1,\ldots k_m}(E)$ and let
$$O_F(b_1,b_2,\ldots) =  \pi_1^*O_{Grass_{k_1}}(b_1)\otimes
\pi_2^*O_{Grass_{k_2}}(b_2)\otimes\ldots$$
Then
$M= O_F(a_1+k_2-k_0, a_2+k_3-k_1, \ldots) $
 is nef or nef and big according to whether
$E$ is.

Let us suppose that $E$ is nef and prove that 
$\F\otimes \BS_+^\lambda(E)\in \ga(X)$ whenever $\F\in \ga_n(X)$ by 
induction.
Suppose that $\F = f_*\omega_{Y}\otimes L$ where $Y$ is smooth and
$L$ is an ample line bundle on $X$. There is no loss of
generality in assuming that $f$ is surjective, since otherwise
we can replace $X$ by its image. 
Consider the Cartesian diagram:

$$\begin{array}{ccc}
F'           &\stackrel{f'}{\longrightarrow} & F\\
p'\downarrow &                               & \downarrow p\\
Y            &\stackrel{f}{\longrightarrow}  & X\\
\end{array}
$$
where  $F' =Flag_{k_1,\ldots k_m}(f^*E)$.
By lemma \ref{lemma:canonical}, 
$$\omega_{F'}\otimes {f'}^*M =
{p'}^*(\omega_Y\otimes det(E)^{\otimes k_m})\otimes O_{F'}(a_1, a_2,\ldots).$$
Therefore 
$${p'}_*(\omega_{F'}\otimes {f'}^*M) = 
\omega_Y\otimes\BS_+^{\lambda}({f}^*E)$$
 
$M\otimes p^*L$ is nef and  also big because the
restriction of $M$ to any fiber is ample so that
$$c_1(M\otimes p^*L)^{dim F} \ge
 c_1(M)^{dim F-dim X}p^*c_1(L)^{dim X} > 0.
$$
Therefore 
$$
f_*\omega_Y\otimes L\otimes \BS_+^\lambda(E)\cong
p_*(f'_*\omega_{F'}\otimes M\otimes p^*L).
$$
is geometrically acyclic by lemma \ref{lemma:ex3}.

Suppose  that we have proved the result
for all sheaves in $\ga_n$. Then clearly it holds
for a direct summand  or an extension  of  sheaves in $\ga_n$.
If $\F = f_*\F'$ with $\F'\in \ga_n$, then $\F\otimes \BS_+^\lambda(E) = 
f_*(\F'\otimes f^*\BS_+^\lambda(E))\in \ga$. 

Now suppose that $E$ is both nef and big. Let $f:Y\to X$ be a
surjective morphism from a smooth projective variety $Y$.
Construct a Cartesian square as above. By lemma \ref{lemma:nef}
$M$ is nef and big, therefore 
$$
f_*\omega_Y\otimes\BS_+^\lambda(E)\cong
p_*(f'_*\omega_{F'}\otimes M)
$$
is geometrically acyclic.
\end{proof}

\begin{remark}
It is possible to improve the theorem slightly by replacing the 
bigness condition on $E$ with the weaker condition that $M$
above is big. This last condition is expressible as the positivity
of an appropriate Chern polynomial of $E$. 
The necessary condition, along with a related vanishing theorem,
can be found in \cite{man1}.
\end{remark}

By a similar argument, we have:

\begin{ex}
Let $f:Y\to X$ be a smooth projective morphism then
$E=f_*\omega_{Y/X}$ is a geometrically semipositive
vector bundle. 
\end{ex}

%%%%%%%%%%%%%%%%
\section{Corollaries of Theorem \ref{thm:positiv}}

We have finished the hard work; now for the entertainment.
The first corollary is really just a reformulation
of the theorem. (See appendix for the definition of $Pos(e)$.)

\begin{cor}
Let $E$ be a rank $e$ vector bundle, and $\lambda\in Pos(e)$. Then
$\BS^\lambda(E)$ is geometrically semipositive if $E$ is nef, and 
$\BS^\lambda(E)$ is geometrically positive if $E$ is nef and big.
\end{cor}

\begin{proof}
Let $\lambda' = (\lambda_1-\lambda_e, \ldots)$ and
$n = \lambda_e-length(\lambda')$, $n$ is positive
by the definition of $Pos(e)$. Then $\BS^\lambda(E) =
\BS_+^\lambda(E)\otimes (det E)^{\otimes n}$. As $det(E)$
is nef, the corollary is immediate.
\end{proof}

\begin{lemma}\label{lemma:highDirIm}
Let $E$ be a geometrically positive vector bundle on a projective
variety $X$. If $f:Y\to X$ is a morphism from a  projective
variety with at worst rational singularities,
then $R^if_*\omega_Y\otimes E$ is geometrically acyclic
for each $i$.
\end{lemma}

\begin{proof}
By the assumption on the singularities, we can replace $Y$ by 
a desingularization without affecting $R^if_*\omega_Y$. Then by
 \cite[2.24]{kollar2}, there exists a
morphism from a smooth projective variety $g:Z\to X$
such that $R^if_*\omega_Y$ is a direct summand
of $g_*\omega_Z$.
\end{proof}

The next result follows directly from the previous 
corollary and lemma.
When the $E_i$ have rank one, this is  essentially 
the Koll\'ar vanishing theorem. For higher rank, this
has some overlap with a theorem
of Manivel \cite{manivel} (which supersedes the
vanishing theorems of \cite{griffiths}, \cite{demailly} and
others - see the introduction to his paper);
the overlap occurs when $Y=X$ is smooth and the $E_i$
are ample.

\begin{cor}\label{cor:van}
If $f:Y\to X$ is a morphism  of projective
varieties such that  $Y$ has at worst
rational singularities and if $E_i$ is a collection
of nef vector bundles and  $\lambda(i)\in Pos(rk(E_i))$  a collection
of partitions, with $i= 1, 2\ldots n$, then
$$H^i(X, R^jf_*\omega_Y\otimes \BS^{\lambda(1)}(E_1)\otimes\ldots
\BS^{\lambda(n)}(E_n)) = 0$$
for all $i>0$ and all $j$, provided that at least one of the
$E_i$ is also big.
\end{cor}

It is possible to define an intermediate notion between
geometric positivity  and semipositivity
similar in spirit to Sommese's  $k$-ampleness
\cite{sommese}. Although we will not develop this systematically
here, we point out a special case:

\begin{lemma}\label{lemma:kample}
Let $f:Y\to X$ be a surjective map of projective varieties
with $Y$ smooth. Let  $E$ be
a geometrically positive vector bundle on $X$, then
$$H^i(Y, \omega_Y\otimes f^*E) = 0$$
for $i > dim\, Y - dim\, X$.
\end{lemma}

\begin{proof}
From the Leray spectral sequence, it suffices to kill
$H^p(R^qf_*\omega_Y\otimes E)$ for all $p+q > k = dim\,X-dim\, Y$.
In fact, these groups vanish for all $p>0$ by the previous 
lemma. And they vanish for $q> k$, by \cite[2.1]{kollar}
\end{proof}

\begin{cor} 
Let $f:Y\to X$ be a surjective map of projective varieties
with $Y$ smooth. Let $E$ be a nef and big
 vector bundle on $X$, then
$$H^i(Y, \omega_Y\otimes \BS^\lambda(f^*E)) = 0$$
for $i > dim\, Y - dim\, X$ and $\lambda\in Pos(rk(E))$.
\end{cor}

Recall that a global section of a vector bundle
is {\em regular}, if its components with respect to any 
local basis forms a regular sequence. The defining
ideal of the scheme of zeros of this section is
locally generated by these components. 

\begin{lemma} Let $X$ be a projective
variety. Let $Z$ be the scheme of zeros of a
regular section of a geometrically positive
vector bundle $E$ on $X$. Let $F^*$ be geometrically
semipositive vector bundle. The restriction
maps 
$$H^i(X, F)\to H^i(Z, F\otimes O_Z)$$
are isomorphisms for $i < dim\, X - rk(E)$
and and injection for $i = dim\, X -rk(E)$.
\end{lemma}

\begin{proof}
Since the section is regular, we have the Koszul
resolution:
$$\ldots \wedge^2E^*\to E^*\to O_X\to O_Z\to 0.$$
This can be broken up into a series of short
exact sequences. $\wedge^kE\otimes F$ is geometrically
positive because it is a summand of $E^{\otimes k}\otimes F$.
Therefore the lemma follows by tensoring these
sequences by $F^*$, and using the vanishing of
$$H^i(X,\wedge^kE^*\otimes F^*)\cong
H^{dimX-i}(X,\omega_X\otimes \wedge^kE\otimes F)^*$$
for $i < dim\, X$
\end{proof}

\begin{cor} With the same notation as in the lemma,
if $X$ is irreducible then 
$Z$ is connected provided that $rk(E) < dim\, X$.
\end{cor}

\begin{proof} $H^0(O_X) \cong H^0(O_Z)$.
\end{proof}
 
\begin{cor}
Let $E$ be a globally generated vector bundle on a smooth
projective variety $X$. Then there is an integer $0\le\kappa\le dim\,X$
such that for any partition $\lambda$, the multiplication map
$$H^0(X, det(E)^{\otimes m})\otimes 
H^i(X, \omega_X\otimes \BS^\lambda(E)\otimes det(E)^{\otimes n})
\to H^i(X, \omega_X\otimes\BS^\lambda(E)\otimes det(E)^{\otimes m+n})$$
is surjective for all $i$, $m\ge 0$, and $n\ge length(\lambda)+\kappa$.
\end{cor}

\begin{proof}
Let $f:X\to G = Grass_r(H^0(E))$ be the map such that $E$
is the pullback of the universal quotient bundle $Q$.
Set $\kappa = dim(f(X))$, and let $i:G\hookrightarrow\PP$
denote the Pl\"ucker embedding.
By the  corollary \ref{cor:van} the Leray spectral sequence collapses
to yield isomorphisms
\begin{equation}\label{eq:hi}
H^i(X,  \omega_X\otimes \BS^\lambda(E)\otimes det(E)^{\otimes n})
\cong H^0(G, R^if_* \omega_X\otimes \BS^\lambda(Q)\otimes O_G(n)) 
\end{equation}
for $n\ge length(\lambda)$. 
Furthermore, the vanishing theorem
implies that 
$$\F = i_*[R^if_* \omega_X\otimes \BS^\lambda(Q)\otimes
O_G(length(\lambda)+\kappa)]$$
is $0$-regular \cite{mumford}.
Therefore by [loc. cit., page 100], 
$$ H^0(\PP, O(m))  \otimes  H^0(\PP, \F) \to  H^0(\PP, \F(m)) $$
is surjective for $m \ge 0$. A simple diagram chase 
using the  maps 
$$ H^0(\PP, O(m))\to H^0(X,det(E)^{\otimes m})$$
and the isomorphism (\ref{eq:hi}) finishes the proof.
\end{proof}

Refinements of the above idea will appear in the forthcoming
work of J.Chipalkatti.

An interesting problem is to find intrinsic criteria for
a vector bundle to be geometrically (semi-)positive. 
Simple examples show that the condition of being
nef and big is not sufficient:

\begin{ex}\label{ex:tangentb} The $i$th  exterior power of the
tangent bundle $E= \wedge^iT_P$ of $P =\PP^n$ is ample.
 However it cannot be geometrically positive (for $i<n$)
because
$$H^{n-i}(P, \omega_P\otimes E) = H^{i}(P, \Omega_P^i)^*\not=0$$
This also shows that $E(-1)$ is not geometrically semipositive
even though it is globally generated (and even ample for $i> 1$). 
\end{ex}

Let us say that a vector bundle $E$ of rank $r$ is {\em strongly
semistable} if and only if  $S^r(E)\otimes det(E)^{-1}$ is nef.
The terminology will be explained shortly. But first let us
point out that one can build simple examples using the following
observations:  $E$ is strongly semistable if and only if $E\otimes L$
is strongly semistable for any line bundle $L$, and a nef vector bundle $E$ is 
strongly semistable if $c_1(E)=0$.

\begin{cor}\label{cor:sstable}
 Let $E$ be a strongly semistable vector bundle.
Then $E$ is geometrically semipositive if $det(E)$ is nef,
and $E$ is geometrically positive if $det(E)$ is nef and big.
\end{cor}

\begin{proof} Since these conditions are
stable by   pullback under a generically finite map,
there is no harm in assuming that the base
variety $X$ is smooth and projective. Let  $r=rk(E)$
then  there exists a smooth variety $Y$  with a line bundle $L$
and a finite map $p:Y\to X$
such that $p^*det(E) = L^{\otimes r}$. The existence of
$Y$ follows from \cite{bloch-gies} or \cite[thm 17]{Ka1}.
 Therefore $S^r(F)$
is nef, where $F = p^*E\otimes L^{-1}$. The Veronese embedding 
$\PP(F)\hookrightarrow \PP(S^r(F))$
shows that  $O_{\PP(F)}(r)$ is also nef. Consequently, so is $F$.
The theorem implies that 
$F = F\otimes det(F) = \BS_+^{(1,0\ldots)}(F)$ is geometrically 
semipositive. Since $E$ is a direct summand of $p_*(F\otimes L)$,
the corollary follows.
\end{proof}

As for the name:

\begin{lemma} If $E$ is a strongly semistable vector
bundle on a smooth projective variety $X$, then
$E$ is semistable with respect to any polarization.
 In other words, given an ample line bundle $H$ and
a torsion free quotient $F$, 
$$\frac{c_1(F)c_1(H)^{dimX-1}}{rk(F)}\ge 
\frac{c_1(E)c_1(H)^{dimX-1}}{rk(E)}.$$
\end{lemma}

\begin{proof}
We will only sketch the proof since the lemma is  not used here.
Proceed as above to construct a  map $p:Y\to X$
of smooth projective varieties such that $p^*det(E)$ has an
$rk(E)$th root $L$. Then $p^*E\otimes L^{-1}$ is nef, hence
also is its restriction to any complete intersection curve $C$
corresponding to a multiple of $p^*H$. The same goes
for $(p^*F\otimes L^{-1})|_C$, so it has nonnegative degree.
This last condition is equivalent to the above inequality.
\end{proof}

Strong semistablitiy is known to be equivalent to ordinary 
semistablity for curves \cite[sect. 3]{miyaoka}. The set
of strongly semistable bundles in an irreducible component
of the moduli space of semistable bundles (for a fixed polarization)
is easily seen to be either empty or the complement of a countable
union of proper subvarieties.

Manivel \cite{man3} has introduced the class of uniformly nef
vector bundles on projective varieties. It is the smallest class 
containing
the bundles  $E\otimes L$, where $E$ is unitary flat and $L$ a nef
line bundle, and which is 
closed under extensions, direct summands, and such that
 $E$ is uniformly nef if and only if its inverse image under a
finite map is.

\begin{cor} A uniformly nef vector bundle is geometrically
semipositive.
\end{cor}

\begin{proof} By  corollary \ref{cor:sstable} and the
comments preceding it,
a vector bundle of the form $E\otimes L$, with $E$ unitary flat
and $L$ a nef line bundle,
 is both nef and strongly semistable and therefore
geometrically semistable. Moreover a vector bundle is clearly
nef and  strongly semistable if and only if its pullback under
a finite map is. This along with the fact that the
 class of geometrically semipositive
vector bundles is closed under extensions and direct summands
implies the corollary.
\end{proof}

Let $X$ be a smooth projective variety.
It is possible to define a countable collection of invariants
for $X$ as follows.
Given a partition $\lambda$ of length less than or equal to $dim\, X$,
define the associated Schur-Weyl invariant by
$$q_\lambda(X) = dim\, H^0(X, \BS^\lambda\Omega_X^1).$$
These numbers, which include the plurigenera,
 are birational invariants. A detailed study  of these invariants
can be found in \cite{man2}.

\begin{cor} Suppose that $\Omega_X^1$ is nef and that $\omega_{X}$
 is big . If
$\lambda$ is a partition such that 
$$\lambda_{dim X}\ge length(\lambda_1-\lambda_{dim X}, \ldots
\lambda_{dim X -1}-\lambda_{dim X}, 0,\ldots ) + 2,$$
then $q_\lambda(X)$  is multiplicative
under \'etale covers and invariant under small deformations.
\end{cor}

\begin{proof}
By our assumptions, there exists a partition $\lambda'\in Pos(dim\,
 X)$,  such that 
$$\BS^\lambda\Omega_X^1 =
\omega_X^{\otimes 2}\otimes \BS^{\lambda'}\Omega_X^1.$$
Consequently, $q_\lambda(X)=\chi(\BS^\lambda\Omega_X^1)$.
The right hand side is 
 multiplicative  and a deformation invariant by Riemann-Roch.
 Furthermore, this equality persists under
\'etale covers because the pullback of a nef (and big) vector bundle
is nef (and big), 
and it persists under small deformations by upper semicontinuity
of cohomology.
\end{proof}

Varieties satisfying the conditions of the previous corollary are
easy to construct by taking subvarieties of an abelian 
 variety of general type.
If $\Omega_{X}^1$ is nef and big (e.g., ample) then the conclusion of the 
corollary holds for any partition satisfying
$$\lambda_{dim X}\ge length(\lambda_1-\lambda_{dim X}, \ldots
\lambda_{dim X -1}-\lambda_{dim X}, 0,\ldots ) + 1.$$

\section{ Blow ups of coherent sheaves, and local algebra}

 Fix a local domain $(R,m)$, with fraction field $K$,
 of essentially
finite type  over $\C$.
When no (or only moderate) confusion is likely,
we  denote an $R$-module and the associated 
quasicoherent sheaf on $spec\, R$ by the same symbol.
 Given a coherent sheaf  $E$ on an integral
scheme $Z$, we will say that
a birational morphism  $f:X\to Z$ {\em frees} $E$
if the quotient of  $f^{*}E$ by its torsion submodule is locally free.

\begin{lemma} Let $Z$ be either an irreducible
 algebraic variety over $\C$
or $spec\, R$.
Given a sequence $E_1,\ldots E_n$ of
coherent  $O_Z$-modules, there exists a resolution
 of singularities $f:X\to Z$ which frees all the $E_i$.
 \end{lemma}
 
\begin{proof} Let $r_i$ be the rank of $E_i$ at the generic
point. 
Let $f_1:X_1\to spec\, R$
be the blow  up of the $r_{1}$st Fitting ideal of $E_1$.
Let $f_{2}:X_2\to spec\, R$  be the blow up of the $r_{2}$nd Fitting ideal of
$f_1^*E_2$, and so on.
By \cite[p. 40]{gruson-raynaud}, $X_n$ frees all of the $E_i$. 
To finish the construction, choose a desingularization $X\to X_n$.
\end{proof}

\setcounter{equation}{0}

\begin{thm} Let $f:X\to spec\, R$ be a resolution of singularities
which frees a sequence $E_1,\ldots E_n$ of finitely generated
$R$-modules, and let $D$ be the exceptional divisor.
Let $\tilde E_i$ be the quotient of $f^*E_i$ by its torsion submodule.
Suppose that $\lambda_i\in Pos(dim\, E_i\otimes K)$. Then
\begin{equation}\label{eq:A}
 H^i(X, \omega_X\otimes \BS^{\lambda_1}(\E_1)\otimes \ldots 
\BS^{\lambda_n}(\E_n)) = 0
\end{equation}
for all $i > 0$, and 
\begin{equation}\label{eq:B}
 H_D^i(X, \BS^{\lambda_1}(\E_1^*)\otimes \ldots 
\BS^{\lambda_n}(\E_n^*)) = 0
\end{equation}
for all $i < dim R$.
\end{thm}

\begin{proof}
We first make a few preliminary observations. (2)
follows from (1) by duality \cite[p 188]{lipman}.
 Next suppose that $\pi:X'\to X$ is a birational map with $X'$
smooth, then the  equality 
${\mathbb R}\pi_*\omega_{X'} = \omega_X$
along with the projection formula shows that it is enough
to prove the (1) result for $X'$.

 The heart of the argument involves globalizing.
 Let $Z$ be a projective
variety with a point $z$ such that $R \cong O_{Z,z}$.
Choose presentations $R^{n_i'}\to R^{n_i}\to E_i\to 0$ for each $i$.
The presentation matrices can be extended to matrices of
regular functions in a neighbourhood of $z$. Hence after blowing
up $Z$ (away from $z$) these matrices extend to maps
$ A_i:O_Z^{n_i'}\to O_Z(D)^{n_i}$ where $D$ is an effective
Cartier divisor containing all the poles of the matrix entries
and such that  $z\notin D$. 
We can assume that $D$ is ample, since otherwise we can replace it with 
 $D+NH$ with $N>> 0$, where $H$ is a very ample divisor
avoiding $z$. Set  ${\mathcal E}_i = coker(A_i)$.
Then by construction ${\mathcal E_i}_z = E_i$.
Let $\bar f: \bar X \to Z$ be a resolution of singularities which frees
the sheaves ${\mathcal E_i}$. By blowing up further, we can
assume that $X' = \bar X\times_Z {spec\, R} $ dominates $X$.
The key point is  that we arranged that each
$\tilde {\mathcal E_i} = \bar f^*{\mathcal E_i}/torsion$
is a  quotient of a direct sum of $f^*O_Z(D)$; as $D$ is ample,
 this implies
that $\tilde {\mathcal E_i}$ is nef and big. Therefore
by theorems \ref{thm:ga} and \ref{thm:positiv},
$$R^i{\bar f}_*(\omega_{\bar X}\otimes \BS^{\lambda_1}(\tilde
{\mathcal E}_1)\otimes \ldots 
\BS^{\lambda_n}(\tilde {\mathcal E}_n)) = 0.$$
We obtain  (1) for $X'$ since cohomology commutes with flat base change.

\end{proof}

\begin{cor}
Suppose that $(R,m)$ is a normal isolated singularity,
and let $E_1,\ldots E_n$ be  a collection of $R$-modules which are
locally free on $spec\,R -\{m\}$. Let $F = 
\BS^{\lambda_1}(E_1)\otimes \ldots 
\BS^{\lambda_n}(E_n)$ where $\lambda_i\in Pos(dim E_i\otimes K)$.
Let $f: X\to spec\, R$ is a resolution of singularities which
frees the $E_i$, and let
$$\tilde F = \BS^{\lambda_1}(f^*E_1/torsion)\otimes \ldots 
\BS^{\lambda_n}(f^*E_n/torsion).$$
Then
$$ H_m^i(F^*) \cong H^{i-1}(X, \tilde F^*)$$
for $2\le i \le dim\, R -1$. In particular,
$$depth\, F^* = 
max\{i\,|\, H^j(X,\tilde F^*) = 0, \forall 1\le j \le i-2\}$$
provided the set on right is nonempty (otherwise $depth\, F^* =2$).
\end{cor}

\begin{proof}     
Consider the spectral sequence
$$ E_2^{pq} = H_m^p(R^qf_*{\tilde F}^*) \Rightarrow H_D^{p+q}(\tilde F^*)$$
If we plot the $E_2$ terms in the $pq$-plane, then our assumptions
imply that the nonzero terms are concentrated along the $p$
and $q$ axes. Since the abutment vanishes for $p+q < dim\, R$,
we can deduce equalities
$$H_m^i(f_*\tilde F^*) = 0$$
for $i < 2$, and isomorphisms
$$H^{i-1}(X, \tilde F^*) \cong E_i^{0,i-1}
\stackrel{\sim}{\to} E_i^{i,0}  \cong H_m^i(f_*\tilde F^*) $$
for $i \ge 2$.
The vanishing of the first two local cohomologies implies that
$f_*\tilde F^*$ is reflexive. Since it coincides with the reflexive
sheaf $F^*$ on the punctured spectrum, it coincides everywhere.
\end{proof}

We now return to the global setting. Let $X$ be a
smooth projective variety.

\begin{lemma}
Let $i:F\to E$ be a map of vector bundles (not necessarily
of constant rank) on $X$ and let $L$ be another line bundle.
Suppose either that $F$ is nef and $L$ is nef and big,
or the other way around.
If $\lambda\in Pos(rk(i(F)\otimes \C(X)))$, then there
exists a geometrically acyclic sheaf $\F$ such that
$$\omega_X\otimes image(\BS^\lambda(F))\otimes L
\subseteq \F\subseteq 
\omega_X\otimes \BS^\lambda(E)\otimes L.$$
The first inclusion is an isomorphism over the open
set where $i(F)$ is locally free.
\end{lemma}

\begin{proof}  
Choose a birational
map $f:Y\to X$ which frees $i(F)$ and with $Y$ smooth. Then
$\tilde F = f^*i(F)/torsion$ is nef or nef and big in accordance
with our assumptions about $F$.
 Therefore
$\F' = \omega_Y\otimes \BS^\lambda(\tilde F) \otimes f^*L$, and
hence also $\F = f_*\F'$, is  geometrically acyclic. There is an 
injection
$\F' \hookrightarrow \omega_Y\otimes \BS^\lambda(f^*E) \otimes f^*L$
 and a surjection
 $\omega_Y \otimes f^*\BS^\lambda(F)\otimes f^*L\to \F'$ which yield the inclusions
 after applying $f_*$.
\end{proof}

\begin{cor} Let $b$ be the dimension of the ``base locus''
 of $E$ i.e. the dimension of the support of $coker[H^0(E)\otimes O_X\to E]$
 (or $0$ if the support is empty).
 If $L$ is a nef and big line bundle and $\lambda\in Pos(rk(E))$,
 then 
 $$H^i(X, \omega_X\otimes \BS^\lambda(E)\otimes L) = 0$$
 for $i > b$.
\end{cor}

%%%%%%%%%%%%%%%%%%

\section{Appendix: Schur-Weyl Functors}
We will briefly describe the Schur-Weyl functors 
\cite[pp 231-7]{fulton-harris} which generalize the
exterior and symmetric powers. A  partition is 
a  nonincreasing sequence of natural numbers 
$\lambda= (\lambda_1, \lambda_2, \ldots)$ which is eventually zero.
Its length  $length(\lambda)$ is the
largest $n$ with $\lambda_n\not= 0$, and $|\lambda| =
\sum_i\lambda_i$. As usual a partition can be represented
by a Young diagram, which consists of $|\lambda|$ boxes arranged
so that $\lambda_i$ boxes lie in the $i$th row.
For each diagram, we will need to choose a labeling of  the boxes by 
integers from $1$ to $|\lambda|$ so as to get a Young tableau; 
for our purposes it will not matter  how this is done. 
Then every partition will
determine an idempotent  $e_\lambda\in
\Q[S_{|\lambda|}]$ obtained by normalizing the so called
 Young symmetrizer. Given a finite
dimensional complex vector space  $E$, define the Schur-Weyl power by
$${\mathbb S}^\lambda E = E^{\otimes |\lambda|}e_\lambda$$
where the symmetric group acts on the right by permuting factors.

 The $\BS^\lambda(E)$ are a complete set of representatives
for  the irreducible representations
of $SL(E)$. According to Borel-Bott-Weil theory,
these representations can be realized geometrically
as spaces of sections of line bundles on the  variety
 $Flag(E)$ of complete flags in  $E$.
Let $\pi_k:Flag(E)\to Grass_k(E)$ be the canonical map to
the Grassmanian of $k$-dimensional {\em quotients} of $E$, and $i_k$ its
Pl\"ucker embedding. Then
$$\BS^\lambda(E) = 
H^0(Flag(E), \bigotimes_k\, \pi_k^*i_k^*O_{\PP(\wedge^kE)}(a_k))$$
where $a_i = \lambda_i-\lambda_{i+1}$. 
This shows, among other things, that $\BS^{\lambda}(E)=0$ as soon
as $length(\lambda)> dim(E)$.  
We can also replace $Flag(E)$, above, by the partial flag variety
$Flag_{k_1, \ldots, k_m}(E)$ parameterizing flags of 
$E\supset E_{k_1}\supset E_{k_2}\ldots $
where $E_{k_i}$ has codimension $k_i$ and $\{k_i\}$ is the set
of the indices where $a_{k_i}\not= 0$. For a treatment along these
lines see \cite[chap. 9]{fulton2}.

All of this makes  perfect sense when $E$ is a replaced
with a vector bundle, provided we replace  $Grass_k$ and
$Flag$ by the corresponding bundles. The validity
of all of the above assertions for $E$ follows by standard base
change arguments.

We will need an extension of the Schur-Weyl functors to 
coherent sheaves where the above formulas are unsuitable.
Define
$$\BS^{\lambda}(E) 
= \frac{E^{\otimes |\lambda|}}{(1-e_{\lambda})E^{\otimes |\lambda|}} 
\cong E^{\otimes |\lambda|}\otimes_{\Q[S_{|\lambda|}]}
\frac{\Q[S_{|\lambda|}]}{(1-e_{\lambda})\Q[S_{|\lambda|}]} .$$
The isomorphism class of the
module $\Q[S_{|\lambda|}]/(1-e_{\lambda})\Q[S_{|\lambda|}]$
is independent of the choice of tableau, so the same goes for 
$\BS^{\lambda}(E)$. This definition
coincides with the previous  one for locally free sheaves,
and takes epimorphisms to epimorphisms, so it can be computed
with the help of a presentation.

It will be convenient to set
$$\BS^{\lambda}_+E = \BS^\lambda E\otimes (det E)^{length(\lambda)}
 = \BS^{\lambda+length(\lambda)\mu}(E)$$
where $\mu = (1,1,\ldots 1, 0\ldots)$ ($rk(E)$ ones).
For each positive integer $e$, let $Pos(e)$ denote the set
of partitions $\lambda$ of length at most $e$, such that
$\lambda_e$ is greater than or equal to the 
length of 
$(\lambda_1-\lambda_e, \ldots, \lambda_{e-1}-\lambda_e, 0,\ldots )$.

%%%%%%%%%%%%%%%%%%%


\begin{thebibliography}{ABC}

\bibitem[BG]{bloch-gies} S. Bloch, D. Gieseker,
{\em The positivity of Chern classes of an ample vector bundle}
Inv. Math. 12 (1971), 112-117
    
\bibitem[D]{demailly} J. P. Demailly, {\em Vanishing theorems for
tensor powers of ample vector bundles}, Inv. Math. 91 (1988), 
203-220

\bibitem[EV]{esnault-viehweg} H. Esnault, E. Viehweg,
{\em Lectures on Vanishing theorems}, Birkhauser (1992)

\bibitem[F]{fulton} W. Fulton, {\em Intersection theory},
Springer-Verlag (1984)

\bibitem[F2]{fulton2} W. Fulton {\em Young tableux},
London Math. Soc. (1997)

\bibitem[FH]{fulton-harris} W. Fulton, J. Harris, {\em Representation theory,
a first course}, Springer-Verlag (1991)

\bibitem[G]{griffiths} P. Griffiths, {\em Hermitean differential geometry,
Chern classes, and positive vector bundles}, Global Analysis, Princeton
(1969)

\bibitem[GR]{gruson-raynaud} L. Gruson, M. Raynaud, {\em Crit\`eres
 de platitude et de projectivit\'e},  Inv. Math. 13
(1971), 1-89

\bibitem[Ha]{hartshorne} R. Hartshorne, {\em Ample vector bundles}
Publ. IHES 29 (1966), 63-94

\bibitem[Hi]{hironaka} H. Hironaka, {Resolution of
singularities of an algebraic variety 
over a field of characteristic 0}, Ann. Math. 79 (1964), 109-326

\bibitem[Ka1]{Ka1} Y Kawamata, {\em Characterization of Abelian varieties,}
Compositio Math. 43 (1981), 253-276


\bibitem[Ka2]{kawamata} Y. Kawamata, {\em A generalization of 
Kodaira-Ramanujam's vanishing theorem}, Math. Ann. 261 (1982) 43-46


\bibitem[Ko1]{kollar} J. Koll\'ar,{\em  Higher direct images of dualizing
sheaves,} Ann. Math. 123 (1986), 11-42

\bibitem[Ko2]{kollar2}J. Koll\'ar, {\em  Higher direct images of dualizing
sheaves II,} Ann. Math. 124 (1986), 171-202

\bibitem[L]{lipman} J. Lipman, {\em Desingularization of two-dimensional
schemes}, Ann. Math. 107 (1978) 151-207

\bibitem[M1]{man1} L. Manivel, {\em Un th\'eor\`eme d'annulation
``\`a la Kawamata-Viehweg''}, Manusc. Math. 83 (1994), 387-404

\bibitem[M2]{man2} L. Manivel, {\em Birational invariants of algebraic
varieties}, J. f. Reine  Angew. Math. 458 (1995) 63-91

\bibitem[M3]{man3} L. Manivel, {\em Th\'eoremes d'annulation sur
 certaines vari\'et\'es projectives}, Comm.  Math. Helv. 71 (1996),
402-425


\bibitem[M4]{manivel} L. Manivel, {\em Vanishing theorems for ample vector
bundles}, Inv. Math. 127 (1997), 401- 416

\bibitem[Mi]{miyaoka} Y. Miyaoka, {\em The Chern classes and
Kodaira dimension of a minimal variety}, Algebraic Geometry, Sendai,
Kinokuniya/North Holland (1987)

\bibitem[Mo]{mori} S. Mori, {\em Classification of higher dimensional
varieties}, Algebraic Geometry, Bowdoin, Amer. Math. Soc. (1987) 

\bibitem[Mu]{mumford} D. Mumford, {\em Curves on an algebraic
surface}, Princeton U. Press (1966)

\bibitem[S]{sommese} A. Sommese, {\em Submanifolds of Abelian 
varieties}, Math. Ann. 233 (1978), 229-256

\bibitem[V]{viehweg} E. Viehweg,{\em  Vanishing theorems}, J. f. Reine
  Angew. Math. 335 (1982), 1-8


    
\end{thebibliography}
\end{document}